\documentclass[a4paper,oneside,11pt]{article}
\usepackage{amssymb}
\usepackage{amsmath}
\usepackage{enumerate}
\usepackage{graphicx}
\usepackage{theorem}
\usepackage{color}
\usepackage{esint}
\usepackage{colortbl}

\usepackage{xcolor}

\usepackage{bbm}

\newtheorem{Theorem}{Theorem}

\newtheorem{Definition}{Definition}

\newtheorem{Lemma}{Lemma}
\newtheorem{Remark}{Remark}

\def\R{{\mathbb R}}

\def\virgp{\raise 2pt\hbox{,}}

\def\XXint#1#2#3{{\setbox0=\hbox{$#1{#2#3}{\int}$ }
\vcenter{\hbox{$#2#3$ }}\kern-.6\wd0}}

\def\<{\langle}
\def\>{\rangle}

\def\({\left(}
\def\){\right)}

\newcommand{\Addressess}{{
  \bigskip
Universidad Aut\'onoma de Madrid, 28049, Spain\\
\textit{Email address:} \texttt{ j.ramos@uam.es}

}}

\makeatletter

\newcommand{\Rmnum}[1]{\expandafter\@slowromancap\romannumeral #1@}
\makeatother

\author{Javier Ramos\thanks{\noindent Research supported by Ram\'on y Cajal grant RYC2018-026354-I (Ministerio de Ciencia, Spain), MICINN grant PID2021-124195NB-C33 and the Severo Ochoa Grant CEX2023-001347-S.}} 
\title{\large \textbf{
A REFINED TRILINEAR KAKEYA ESTIMATE IN $\mathbb{R}^3$}}
\date{}
\begin{document}
\maketitle
\begin{abstract}
We prove a refined trilinear Kakeya estimate in three dimensions, valid for small values of the transversality parameter.\end{abstract}

\section{Introduction}
The multilinear Kakeya problem was first introduced in \cite{ben}. The problem can be stated as follows: given $d$ families $\mathcal{T}=\{\mathcal{T}_n\}_{n=1}^d$ of tubes with unit width, such that for each tube $T_n\in\mathcal{T}_n$, the direction $d(T_n)\in\mathbb{S}^{d-1}$ lies within a small fixed neighborhood of the $n$’th standard basis vector $e_n$, do we have the inequality 
\begin{align}\label{mainK}
\int_{\R^d} \prod_{n=1}^d \Big( \sum_{T_{n}\in \mathcal{T}_n } \chi_{T_{n}}\Big)^{q}\leq C_{q,n} \prod_{n=1}^d |\mathcal{T}_n|^{q}
\end{align}
for $q\geq \frac{1}{d-1}$?. 
The conjecture, up to the endpoint, was proven in \cite{bct} using a heat-flow technique. The endpoint case $q=\frac{1}{d-1}$ was resolved by Guth \cite{G} using methods from algebraic topology (see also \cite{cv} for an alternative proof based on the Borsuk-Ulam theorem). This endpoint result included a refinement: it permits the families $\{\mathcal{T}_n\}_{n=1}^d$ to satisfy a more general transversality condition. Specifically, given $d$ families of tubes $\{\mathcal{T}_n\}_{n=1}^d$ of unit width such that for each tuple $(T_1, \cdots, T_d)\in\mathcal{T}_1\times \cdots\times \mathcal{T}_d$ we have $$|d(T_1)\wedge \cdots \wedge d(T_d))|\sim \theta,$$ then the following inequality holds
\begin{align}\label{GK1}
\int_{\R^d} \prod_{n=1}^d \Big( \sum_{T_{n}\in \mathcal{T}_n } \chi_{T_{n}}\Big)^{\frac{1}{d-1}}\leq\theta^{-\frac{1}{2}} \prod_{n=1}^d |\mathcal{T}_n|^{\frac{1}{d-1}}.
\end{align}
In \cite{gu2}, a simplified proof of the multilinear Kakeya inequality \eqref{mainK} was given, using a clever induction-on-scales argument. However, this approach neither establishes the endpoint case nor yields the sharp dependence on the transversality parameter.

One of the principal motivations for studying multilinear Kakeya estimates lies in their crucial role in linear problems. In particular, the groundbreaking broad-narrow decomposition technique introduced by Bourgain--Guth \cite{boguth} provides a powerful framework for exploiting these multilinear estimates to derive linear estimates; see, for example, \cite{bode}, or \cite{sticky}, which led to the recent breakthrough in the proof of the Kakeya set conjecture in three dimensions \cite{conje}. See also, for example, \cite{borg:cone}, \cite{cordoba2} or \cite{tvv} for the bilinear-to-linear approach.

There are other generalizations, including replacing neighborhoods of lines with higher-dimensional affine subspaces (see, for example, \cite{zhang} or \cite{Rog}), or considering levels of multilinearity that do not match the ambient dimension (see, for example, \cite{boguth}, \cite{guza}, \cite{hrz} or \cite{za}).

The main contribution of this paper is a refinement of the multilinear Kakeya estimate in three dimensions. Informally, before introducing the necessary notation, it states that under the hypotheses of the estimate \eqref{GK1}, and denoting $Q_R$ as the ball of radius $R$, we have:
\begin{align}\label{me}
\int_{Q_R \cap X} \prod_{n=1}^3 \Big( \sum_{T_{n}\in \mathcal{T}_n } \chi_{T_{n}}\Big)^{\frac{1}{2}}\leq  F(X, \mathcal{T}) \prod_{n=1}^3 |\mathcal{T}_n|^{\frac{1}{2}},
\end{align}
where $X$ is a union of unit balls, and the quantity $F(X, \mathcal{T})$ satisfies $F(X, \mathcal{T}) \leq \theta^{-\frac{1}{2}}$. The value of $F(X, \mathcal{T})$ depends on how the set $X$ is organized into subsets determined by the transversality at each scale. In particular, there are two distinguished scales: the initial one, given by a parallelepiped $P$ whose side lengths and orientation depend on $\theta$, and the final one, which is the ball of radius $R$. If $X$ has density  $\ll 1$ inside $P$, then $ F(X, \mathcal{T}) \ll \theta^{-\frac{1}{2}}$. Similarly, if $X$ has density $\sim 1$ inside $Q_R$, we also obtain $ F(X, \mathcal{T}) \ll \theta^{-\frac{1}{2}}$. The contribution from intermediate scales is more subtle and requires comparing pairs of scales. Additionally, at each scale, $F(X, \mathcal{T})$ accounts for the tubes that do not intersect $X$, in order to extract further gain.

Some other refinements of multilinear and linear estimates have been established in the context of Kakeya and restriction theory; see, for example, \cite{du1}, \cite{du2} or \cite{li}. However, it appears that this is the first instance in which the refinement is governed by the transversality, in order to fully exploit the potential of the multilinear estimates in the small transversality regime.

The proof of our refinement \eqref{me} will make use of Guth trilinear estimate \eqref{GK1}; however, in the final section, we will present an alternative, simpler proof of a weaker version of \eqref{GK1}, which includes an additional factor $C_\epsilon R^{ \epsilon} $. In particular, at the cost of this factor, we will obtain \eqref{me} via a self-contained proof.

The equivalence between the multilinear restriction and the multilinear Kakeya problem directly yields a similar refinement for the multilinear restriction. However, the standard argument incurs a $\log R$ loss, which we will avoid for sufficiently small transversality in \cite{ra2}.

\section{Notation and definitions}

We begin by introducing the subsets of directions in which the families of tubes  $\{\mathcal{T}_n\}_n$ point. These definitions are the same as those given in \cite{ra}, in the context of a trilinear restriction estimate.

\begin{Definition}
Let $S_1,S_2,S_3\subset  \{(\xi,\tfrac{1}{2}|\xi|^2):\;|\xi|\leq 2\}$ satisfy \begin{align}\label{trans}|n(\xi_1)\wedge n(\xi_2) \wedge n(\xi_3)|\sim \theta\end{align} for all choices $\xi_n\in S_n$, where $n(\xi_n)$ is the normal vector to $S_n$ in $\xi_n$.  Let the collection of tubes $\{\mathcal{T}_n\}_n$ be such that for every $T_n\in\mathcal{T}_n$, we have $d(T_n)=n(\xi_n)$ for some $\xi_n\in S_n$. We say that $\{\mathcal{T}_n\}_n$ is of $(r,j,t,w,m)$-type if there exist $(r,j,t,w,m)\in \mathbb{N}\times \mathbb{N}\times \mathbb{N}\times \mathbb{S}^1\times \mathbb{R}^2$, $r\leq j$, such that, possibly after reordering the sets $S_n$, the following holds:
\begin{align*}
&S_1\subset \{(\xi,\tfrac{1}{2}|\xi|^2):\;\xi\in\tau^{j}_k \cap \mathfrak{t}^{j+t}_{w,m}\},\\
&S_2\subset \{(\xi,\tfrac{1}{2}|\xi|^2):\;\xi\in\tau^{j}_{k'} \cap \mathfrak{t}^{j+t}_{w,m}\},\\
&S_3\subset\{(\xi,\tfrac{1}{2}|\xi|^2):\;\xi\in\tau^{r}_{k''} \cap \mathfrak{t}^{r+t}_{w',m}\},
\end{align*}
for some $k,k',k'',w'$ such that $m\in \tau^{j}_k$, $  d(\tau^{j}_k,\tau^{j}_{k'})\sim 2^{-j}$, $ d(\tau^{j}_{k'},\tau^{r}_{k''})\sim d(\tau^{j}_k,\tau^{r}_{k''})\sim 2^{-r}$, $ |w-w'| \sim 2^{-t}$ and 
\begin{align*}
\theta\sim 2^{-j} 2^{-r} 2^{-t}.
\end{align*}
\end{Definition}

Here, $\tau^j_{k}$ denotes the square of side length $2^{-j}$ whose lower-left vertex is located at the point $k$, and $\mathfrak{t}^{j}_{w,m}$ denotes the strip in the plane of width $2^{-j}$, oriented in the direction $w$, and passing through the point $m\in [0,2]^2$. That is,
$$\mathfrak{t}^{j}_{w,m}:=\big\{\xi\in [0,2]^2: \big|\xi-m-w\,((\xi-m)\cdot w)\big|\leq 2^{-j}\big\}.$$

In \cite{ra}, it was shown how transversal sets of directions satisfying \eqref{trans} can be reduced to the case of some $(r,j,t,w,m)$-type.

We now introduce the parallelepipeds that we will be working with, determined by the parameters $(j,t,w,m)$. These parallelepipeds are the appropriate ones for carrying out the induction; see \cite{ra} for further details.  
We note that, in contrast to those defined in \cite{ra}, the parallelepipeds considered here are rescaled by a factor of $2^{-(j+t)}$.

\begin{Definition}
Let $(j,t,w,m,E,F)\in \mathbb{N}\times \mathbb{N}\times \mathbb{S}^1\times \mathbb{R}^2\times (0,1]\times (0,1]$. Define a parallelepiped with spanning vectors $e_1=(w,0)$, $e_2=\big((w,0)\times(m,-1)\big)$, $e_3=(m,-1)$ and correponding side lengths $s_1=Ej+t$, $s_2=Ej+Ft$, $s_3=j+t$. We denote this parallelepiped by $\mathcal{P}(j,t,w,m,E,F)$. The rescaled parallelepiped $ 2^{\lambda} \mathcal{P}(j,t,w,m,E,F)$ is denoted by $\mathcal{P}(j,t,w,m,E,F)[\lambda]$. We also write $\mathcal{P}(j,t,w,m,E,F)[\lambda](a)$ for the translation of this parallelepiped centered at the point $a\in \R^3$.
\end{Definition}

The case $\mathfrak{P}(j,t,w,m,E,F)[Mt+Nj]$ will be denoted by $\mathfrak{P}(j,t,E,F)[M,N]$, and when $E=F=0$, we will simply write $\mathfrak{P}(j,t)[M,N]$.


As part of the argument, we will move from one scale to another, partitioning parallelepipeds into smaller ones. The following definition will be useful for this purpose:
\begin{Definition}
Let $B,C$ be parallelepipeds in $\R^3$.  We define $\mathcal{A}(B,C)$ as a collection of disjoint translates of $C$ covering $B$, that is, a subset $\{C+a:\; a\in\R^3\}$, satisfying:
$$\bigcup_{C'\in \mathcal{A}(B,C)} C'\cap B=B,$$
 $C_1\cap C_2=\emptyset$ for all distinct $C_1,C_2\in \mathcal{A}(B,C)$, and there exists $C'\in \mathcal{A}(B,C)$ whose center coincides with the center of $B$.
\end{Definition}



The collection of unit-width tubes in all possible directions and positions will be denoted by $\mathfrak{T}$. For any $\mathcal{T} \subset \mathfrak{T}$, we denote by $\mathcal{T}[\lambda]$ the same family of tubes $\mathcal{T}$, but with radius $2^{\lambda}$ instead of $1$.
 We now introduce the following subsets of the family $\mathfrak{T}$:
\begin{Definition}
Given subsets $C,B\subseteq Q_R$, $A\subseteq \mathcal{A}(Q_R,C)$, and a collection $\mathcal{T}_n \subseteq \mathfrak{T}$, we define:
\begin{align*}&\mathcal{T}_n(C):=\{ T_n\in \mathcal{T}_n:\quad T_n\cap C\neq \emptyset \},\\
&\mathcal{T}(T_{n}, C):=\big\{T_{n}'\in \mathfrak{T}(C), \; d(T_n')=d(T_n)\big\},\\
&\mathcal{T}_n(A, B):=\{T_n\in \mathcal{T}_{n}(T_{n}',C_i):\; T_n'\in \mathcal{T}_{n}(C_i), \; C_i\in A\cap B\}.
\end{align*}
\end{Definition}

Note that for a tube $T_n \in \mathcal{T}_n$, the set $\mathcal{T}(T_{n}, C)$ may include tubes $T' \in \mathfrak{T}(C)$ with $d(T') = d(T_n)$ but $T' \notin \mathcal{T}_n$. Additionally, the union $\bigcup_{T_n \in \mathcal{T}_n(C)} \mathcal{T}(T_n, C)$ may have multiplicity $>1$, as distinct tubes $T_n, T_n' \in \mathcal{T}_n(C)$ with $d(T_n) = d(T_n')$ satisfy $\mathcal{T}(T_n, C) = \mathcal{T}(T_n', C)$.

In order to increase the scale, we introduce the following parameters:  $\pmb{M}=\{M_1,\cdots, M_{S-1}\}$, $\pmb{N}=\{N_1,\cdots, N_{S-1}\}$  where $S$ is chosen such that $M_{S-1} (j+1)+N_{S-1} (t+1)=\log R$, $M_1=N_1=0$, and for some $\epsilon>0$, $\max(\epsilon,\frac{r}{j})\leq M_{s+1}-M_s\leq 1$ and $\epsilon\leq N_{s+1}-N_s\leq 1$ for $s=1,\cdots, S-2$.
For $2\leq s\leq S-1$, we define the increments $E_s=M_{s}-M_{s-1}$ and $F_s=N_{s}-N_{s-1}$, with initial values  and $F_1=E_1=E_{S}=F_{S}=1$. 

We will use the following variables to denote the density of the parallelepipeds at each scale: $\pmb{\mu}=\{\pmb{\mu}_s\}_{s=1}^{S}=\{(\mu_{1,s},\mu_{2,s},\mu_{3,s},\mu_{4,s},\mu_{5,s})\}_{s=1}^{S}.$ 
To quantify the proportion of tubes that do not intersect a parallelepipeds at one scale but do intersect the parallelepiped at the next scale, we introduce the following variables: $\pmb{\beta}=\{\pmb{\beta}_{s}\}_{s=1}^{S}=\{\{(\beta_{n,1,s},\beta_{n,2,s})\}_{n=1}^3\}_{s=1}^{S}$. We will use the notation $\overline{\pmb{\beta}}_{i,s}=\prod_{n=1}^3 \beta_{n,i,s}$, for $i=1,2$.



\begin{Definition}
Let $(R_1, R_2, R_3)$ and $C$ be a set in $\R^3$, we define \begin{align*}\mathfrak{A}[(R_1, R_2, R_3),\mathcal{T}, C]=\bigcap_{n=1}^3 \mathfrak{A}[R_n,\mathcal{T}_n,C],\end{align*} where
\begin{align*}
\mathfrak{A}[R_n,\mathcal{T}_n, C]:= \{C\in  \mathcal{A}(Q_R,C):\;  |\mathcal{T}_n(C)|\sim R_{n}\}.\end{align*}
Let $\textbf{R}=\{\textbf{R}_{1,s}, \textbf{R}_{2,s}\}_{s=1}^S=\{\{R_{n,1,s}\}_{n=1}^3, \{R_{n,2,s}\}_{n=1}^3\}_{s=1}^S$,  be some dyadic numbers, $\mathcal{T}=\{\mathcal{T}_n\}_{n=1}^3 \subseteq \mathfrak{T}$ and
\begin{align*}\mathcal{T}_{n,s+1}&=\mathcal{T}_n[M_{s+1}j+N_{s+1}t]( \mathcal{A}(Q_R,\mathfrak{P}(j,t,E_{s+1},F_{s+1})[M_s, N_s]), Q_R).
\end{align*}
We denote $\mathfrak{A}[1]=\mathfrak{A}[\textbf{R}_{1,1},\mathcal{T}, Q_1]$, $\mathfrak{P}[1]=\mathfrak{A}[\textbf{R}_{2,1},\mathcal{T}, Q_1]$ and for $1\leq s\leq S-1$,\begin{align*}&\mathfrak{A}[s+1]=\mathfrak{A}[\textbf{R}_{1,s+1},\mathcal{T}_{s+1}, \mathfrak{P}(j,t,E_{s+1},F_{s+1})[M_s, N_s]],\\&\mathfrak{P}[s]=\mathfrak{A}[\textbf{R}_{2,s},\mathcal{T}_{s},\mathfrak{P}(j,t)[M_{s}, N_{s}]], \; s\neq 1.\end{align*}
\end{Definition}

By the geometry of our parallelepipeds, a subset $B\subset\mathcal{A}( \mathfrak{P}(j,t,E_{s+1},0)[M_s,N_{s}](a), \mathfrak{P}(j,t)[M_s,N_{s}])$ can be intersected by a number of tubes in $\mathcal{T}_3$, ranging from $R_{3,2,s}$ up to $R_{3,2,s} 2^{E_{s+1} j+r}$. We quantify this observation with the following definition:

\begin{Definition}
For subset $B\subset  \mathcal{A}( \mathfrak{P}(j,t,E_{s+1},0)[M_s,N_{s}](a), \mathfrak{P}(j,t)[M_s,N_{s}])\cap \mathfrak{P}[s]$, we define $1\leq L(B)\leq 2^{E_{s+1}j+r}$, as the quantity satisfying $|\mathcal{T}_3(B)|\sim L(B)R_{3,2,s}$.

\end{Definition}

We can now define the subset of $Q_R$ where we integrate to obtain the refinement:
\begin{Definition}
Fix $r, j, t, w, m, \pmb{R}, \pmb{M}, \pmb{N}, \pmb{\mu}, \pmb{\beta}$. Let $1\leq s_0\leq S$, and let $\mathfrak{A}^*=\{\mathfrak{A}^*[s]\}_{s=1}^{S}$, $\mathfrak{P}^*=\{\mathfrak{P}^*[s]\}_{s=1}^{S-1}$ be such that $\mathfrak{A}^*[s]\subseteq\mathfrak{A}[s]$ and $\mathfrak{P}^*[s]\subseteq\mathfrak{P}[s]$. Let $\{\mathcal{T}_n\}_{n=1}^{3}$ be a collection of tubes of $(r,j,t,w,m)$-type.  We define
\begin{align*}
\mathfrak{B}_{s_0}:= \bigcap_{s=s_0}^{S}\mathfrak{A}_{1,s}\cap \mathfrak{A}_{2,s},
\end{align*}
where the sets $\big\{\mathfrak{A}_{1,s},\mathfrak{A}_{2,s}\big\}_{s=1}^{S}$ are defined as follows:

Set $\mathfrak{A}_{2,1}=\mathfrak{A}^*[1]$ and  $\mathfrak{A}_{1,S}=Q_R$. For  $1\leq s\leq S-1$, define $$\mathfrak{A}_{1,s}=  \mathfrak{A}^1_{1,s} \cap \mathfrak{A}^2_{1,s} \cap \mathfrak{P}^*[s],$$ where 
 \begin{align*}
 \mathfrak{A}^1_{1,s} :=\Big\{&P\in \mathcal{A}(Q_R,\mathfrak{P}(j,t)[M_s, N_s]):\;\; \\
 &\#\Big\{P_1\in  \mathcal{A}(P,\mathfrak{P}(j,0)[M_s, N_s ]):\;\; \\
&\#\big\{P_2\in  \mathcal{A}(P_1,\mathfrak{P}(r,0)[M_s,N_s ]): \; \\
&\#\{P'\in \mathcal{A}(P_2,\mathfrak{P}(0,0)[M_s,N_s ]):\; \mathfrak{A}_{2,s}\cap P'\neq \emptyset\} \sim \mu_{3,s}\big\} \sim 2^{j(1-E_{s})}
 \mu_{2,s}\Big\}\sim 2^{2t(1-F_{s})}\mu_{1,s}\Big\},
 \end{align*}
 and
 \begin{align*}
 \mathfrak{A}^2_{1,s}:=\Big\{P\in \mathcal{A}(Q_R,\mathfrak{P}(j,t)[M_s, N_s]):\; \frac{\#\{T_{i}\in\mathcal{T}_n:\quad T_{i} \cap P  \cap\mathfrak{A}_{2,s}\neq \emptyset\}}{\#\{T_{i}\in \mathcal{T}_n:\quad T_{i} \cap P\neq \emptyset\}}\sim \beta_{n,1,s}\quad n=1,2,3\Big\}.
\end{align*}
For $1\leq s\leq S-1$, define $$\mathfrak{A}_{2,s+1}=\mathfrak{A}^1_{2,s+1}\cap \mathfrak{A}^2_{2,s+1}\cap \mathfrak{A}^*[s+1],$$ where 
 \begin{align*}
\mathfrak{A}^1_{2,s+1}:=\Big\{&P_3\in \mathcal{A}(Q_R,\mathfrak{P}(j,t,E_{s+1},F_{s+1})[M_s, N_{s}]):\; \\
&\#\big\{ P_4\in  \mathcal{A}(P_3,\mathfrak{P}(j,t,E_{s+1},0)[M_s, N_{s}]):\quad
& L(P_4\cap \mathfrak{A}_{1,s})\sim\mu_{5,s+1}\big\}\sim \mu_{4,s+1}\Big\},
\end{align*}
and
\begin{align*}
 \mathfrak{A}^2_{2,s+1}:=\Big\{P\in \mathcal{A}(Q_R,&\,\mathfrak{P}(j,t,E_{s+1},F_{s+1})[M_s, N_s]):\; \\
 &\frac{\#\{T_{i}\in \mathcal{T}_n:\quad T_{i} \cap P  \cap \mathfrak{A}_{1,s}\neq \emptyset\}}{\#\{T_{i}\in \mathcal{T}_n:\quad T_{i} \cap P\neq \emptyset\}}\sim\beta_{n,2,s+1}\quad n=1,2,3\Big\} .
\end{align*}

\end{Definition}

\begin{Remark}
We observe that in the previous definition we just considered $(\pmb{\mu}, \pmb{\beta})$ such that  for every $1\leq s \leq S$, $n\in\{1,2,3\}$, 
\begin{align*}
&\mu_{1,s}\leq 2^{2 F_s t},\\
&\mu_{2,s}\leq 2^{E_s j-r},\\  
&\mu_{3,s}\leq 2^{r},\\  
&\mu_{4,s}\leq 2^{F_s t},\\ 
&\mu_{5,s}\leq 2^{E_s j+r},\\
&\beta_{n,1,s},\beta_{n,2,s}\leq 1.\end{align*}
In particular, we have
\begin{align*}
&\mu_{1,1}^{\frac{1}{4}}\mu_{2,1}^{\frac{1}{2}}\mu_{3,1}\leq 2^{\frac{1}{2}(j+r+t)}=\theta^{-\frac{1}{2}}.\end{align*}

\end{Remark}



The refinement factors are given by
\begin{align*}
\mathcal{S}(\pmb{\beta},\pmb{\mu},s)&=  \min(1,\overline{\pmb{\beta}}_{1,s}^{\frac{1}{2}}  \overline{\pmb{\beta}}_{2,s}^{\frac{1}{2}} \mu_{1,s}^{\frac{1}{4}} \mu_{2,s}^{\frac{1}{2}}\mu_{3,s}\mu_{4,s}^{-\frac{1}{2}}\mu_{5,s}^{-\frac{1}{2}})\quad 2\leq s\leq S-1,\\
\mathcal{S}(\pmb{\beta},\pmb{\mu},1)&=  \overline{\pmb{\beta}}_{1,1}^{\frac{1}{2}}   \mu_{1,1}^{\frac{1}{4}} \mu_{2,1}^{\frac{1}{2}}\mu_{3,1},\\
\mathcal{S}(\pmb{\beta},\pmb{\mu},S)&=  \overline{\pmb{\beta}}_{2,S}^{\frac{1}{2}}\mu_{4,S}^{-\frac{1}{2}}\mu_{5,S}^{-\frac{1}{2}}.
\end{align*}


The last definition is for the induction-on-scale argument:
\begin{Definition}
Fix $r, j, t, w, m, \pmb{R}, \mathfrak{A}^*,\mathfrak{P}^*,\pmb{M}, \pmb{N}, \pmb{\mu}, \pmb{\beta}$. Let $\{\mathcal{T}_n\}_n$ be a collection of tubes of $(r,j,t,w,m) $-type. For every $1\leq s\leq s'\leq S$, we denote by $\mathcal{K}(s,s')$ the smallest constant $C$ such that 
\begin{align*}
 \sum_{B_{s}\in  \mathfrak{B}_{s}}\int_{P_{s'}\cap B_{s}} \prod_{n=1}^3 \Big(\ \sum_{T_{n}\in \mathcal{T}_n[M_{s}j+N_{s}t]( \mathfrak{B}_{s}, P_{s'}) } \chi_{T_{n}}\Big)^{\frac{1}{2}}\leq C 2^{3(M_{s}j+N_{s}t)}  \prod_{n=1}^3 |\mathcal{T}_n[M_{s}j+N_{s}t]( \mathfrak{B}_{s},P_{s'})|^{\frac{1}{2}},
\end{align*}
for every $P_{s'}\in  \mathfrak{A}_{1,s'}$. 
\end{Definition}

Note that for the case $s=1$, the inequality simplifies to
\begin{align*}
\sum_{B_{1}\in  \mathfrak{B}_{1}}\int_{P_{s'}\cap B_{1}} \prod_{n=1}^3 \Big(\ \sum_{T_{n}\in \mathcal{T}_n } \chi_{T_{n}}\Big)^{\frac{1}{2}} \leq C  \prod_{n=1}^3 |\mathcal{T}_n(P_{s'})|^{\frac{1}{2}}.\end{align*}
\subsection{Statement of the main result}

We set
\begin{align*}
I= \{s:\; 1\leq s\leq S-1,\;\; \mathcal{S}(\pmb{\beta},\pmb{\mu},s+1)\leq C^{-1}\}\cup\{S\},
\end{align*}
where $C$ is a constant. We order the elements of $I$ and write $I=\{s_1,s_2,\cdots,s_{|I|}\}$, where $s_{|I|}=S$. 
We are now ready to state the refinement of the trilinear Kakeya estimate:
 \begin{Theorem}\label{GTeM}
 $$ \mathcal{K}(1,S)\leq C^{|I|}\mathcal{S}(\pmb{\beta},\pmb{\mu},1)\prod_{s_i\in I\setminus S}\mathcal{S}(\pmb{\beta},\pmb{\mu},s_i+1). $$
 \end{Theorem}
It follows that:
\begin{align*}
\sum_{B_{1}\in  \mathfrak{B}_{1}}&\int_{B_1\cap Q_R} \prod_{n=1}^3 \Big( \sum_{T_{n}\in \mathcal{T}_n } \chi_{T_{n}}\Big)^{\frac{1}{2}}\leq C^{|I|}\mathcal{S}(\pmb{\beta},\pmb{\mu},1)\prod_{s_i\in I\setminus S}\mathcal{S}(\pmb{\beta},\pmb{\mu},s_i+1)\prod_{n=1}^3 |\mathcal{T}_n|^{\frac{1}{2}}.
\end{align*}
\begin{Remark}
We observe that if $|I|=1$, then the right-hand side in Theorem \ref{GTeM} is simply $C \mathcal{S}(\pmb{\beta},\pmb{\mu},1)$.
Moreover, we always have  $$C^{|I|}\mathcal{S}(\pmb{\beta},\pmb{\mu},1)\prod_{s_i\in I\setminus S}\mathcal{S}(\pmb{\beta},\pmb{\mu},s_i+1)\lesssim  \mathcal{S}(\pmb{\beta},\pmb{\mu},1)\mathcal{S}(\pmb{\beta},\pmb{\mu},S)\leq \theta^{-\frac{1}{2}}.$$
\end{Remark}

We now extend the discussion from the introduction regarding the philosophy behind the refinement. The quantity $\mathcal{S}(\pmb{\beta},\pmb{\mu},1)$ corresponds to the refinement with respect to $\theta^{-\frac{1}{2}}$. When $\mathfrak{B}_1$ has full density inside $\mathfrak{P}(j,t)$, we have $\mathcal{S}(\pmb{\beta},\pmb{\mu},1)=\theta^{-\frac{1}{2}}$; however, when the density is lower, we obtain a corresponding gain $\mathcal{S}(\pmb{\beta},\pmb{\mu},1)= \big(\theta^{\frac{1}{2}} \mu_{1,1}^{\frac{1}{4}} \mu_{2,1}^{\frac{1}{2}}\mu_{3,1}\big) \theta^{-\frac{1}{2}}$.  The factor $\mathcal{S}(\pmb{\beta},\pmb{\mu},S)$ satisfies $\mathcal{S}(\pmb{\beta},\pmb{\mu},S)\sim 1$  when $|\mathfrak{B}_{S-1}\cap Q_R|\sim 1$ and becomes smaller as the density increases, again yielding a gain $\mathcal{S}(\pmb{\beta},\pmb{\mu},S)\sim\mu_{4,s}^{-\frac{1}{2}}\mu_{5,s}^{-\frac{1}{2}}$. For the intermediate scales, we must combine information from two consecutive scales, which may yield a gain for certain configurations. 

The parameters $\pmb{\beta}$ yield a gain at the first scale when the tubes $\mathcal{T}$ intersect $\mathfrak{A}_{1,1}$ but not $\mathfrak{A}_{2,1}$, and yield a gain at the final scale when the tubes $\mathcal{T}$ do not intersect $\mathfrak{A}_{1,S-1}$. For the intermediate scales the gain will depende on consecutive scales and on $\pmb{\mu}$. In particular, we notice that it is not possible to replace $\mathcal{S}(\pmb{\beta},\pmb{\mu},s)$ by
\begin{align*}
\overline{\pmb{\beta}}_{1,s}^{\frac{1}{2}}  \overline{\pmb{\beta}}_{2,s}^{\frac{1}{2}}\min(1, \mu_{1,s}^{\frac{1}{4}} \mu_{2,s}^{\frac{1}{2}}\mu_{3,s}\mu_{4,s}^{-\frac{1}{2}}\mu_{5,s}^{-\frac{1}{2}})\quad 2\leq s\leq S-1.\end{align*}   
For any $1\leq s\leq S$, the absence of gain in $\mathcal{S}(\pmb{\beta},\pmb{\mu},s)$ provides structural information about the distribution of the tubes.


\subsection{Proof of Theorem \ref{GTeM}.}

We will prove Theorem~\ref{GTeM} using an induction argument.

\begin{Lemma}\label{L1}
For all $1\leq s< s'\leq S$, and also for $1\leq s= s'\leq S-1$, the following holds:
 $$ \mathcal{K}(s,s')\lesssim  \overline{\pmb{\beta}}_{1,s}^{\frac{1}{2}} \mu_{1,s}^{\frac{1}{4}} \mu_{2,s}^{\frac{1}{2}} \mu_{3,s}\big(2^{(1-F_{s})t}\big)^{\frac{1}{2}} \big(2^{(1-E_{s})j}\big)^{\frac{1}{2}}.$$ 

\end{Lemma}
\begin{Lemma}\label{L2}
If $|I|> 1$, for every $s_i\in I$, $1\leq i\leq |I|-1$,
\begin{align*}\mathcal{K}(1,s_{i+1})\lesssim  \mathcal{K}(1,s_i) \mathcal{S}(\pmb{\beta},\pmb{\mu},s_{i}+1).
\end{align*}
\end{Lemma}

\textbf{Proof of Theorem \ref{GTeM}.}
We observe that Lemma \ref{L1} implies that  $\mathcal{K}(1,s)\lesssim  \mathcal{S}(\pmb{\beta},\pmb{\mu},1)$ for all $1\leq s \leq S$. In particular, $\mathcal{K}(1,s_1)\lesssim  \mathcal{S}(\pmb{\beta},\pmb{\mu},1)$. If $|I|=1$ we are done. Otherwise, applying Lemma \ref{L2} successively $|I|-1$ times, we obtain the desired result:
\begin{align*}\mathcal{K}(1,S)&\leq C \mathcal{K}(1,s_{|I|-1}) \mathcal{S}(\pmb{\beta},\pmb{\mu},s_{|I|-1}+1)\\
&\leq \cdots \\
&\leq C^{|I|-1}  \mathcal{K}(1,s_1)\prod_{1\leq i\leq  |I|-1}\mathcal{S}(\pmb{\beta},\pmb{\mu},s_i+1)\\
&\leq C^{|I|} \mathcal{S}(\pmb{\beta},\pmb{\mu},1)\prod_{1\leq i\leq  |I|-1}\mathcal{S}(\pmb{\beta},\pmb{\mu},s_i+1).
\end{align*}

\begin{flushright}
$\blacksquare$
\end{flushright}

\textbf{Proof of Lemma \ref{L2}.}



We denote for $B_{s_i+1}\in\mathfrak{A}_{2,s_i+1}$,
$$N(B_{s_i+1})=\#\{P\in \mathfrak{A}_{1,s_i}\cap B_{s_i+1}\}.$$
Now, using the definitions of $N(B_{s_i+1})$ and $\mathcal{K}(s_i)$, we have for $P_{s_i}\in  \mathfrak{A}_{1,s_i}\cap B_{s_i+1}$,
\begin{align*}
 \sum_{B_{1}\in \mathfrak{B}_1} \int_{  B_{s_i+1}\cap B_{1}} \prod_{n=1}^3 \Big( \sum_{T_{n}\in \mathcal{T}_n } \chi_{T_{n}}\Big)^{\frac{1}{2}}&\sim N(B_{s_i+1}) \sum_{B_1\in \mathfrak{B}_1} \int_{B_{s_i+1}\cap P_{s_i} \cap B_{1}} \prod_{n=1}^3 \Big( \sum_{T_{n}\in \mathcal{T}_n} \chi_{T_{n}}\Big)^{\frac{1}{2}}\\
&\lesssim N(B_{s_i+1}) \mathcal{K}(1,s_i) \prod_{n=1}^3 |\mathcal{T}_n(P_{s_i})|^{\frac{1}{2}}.
\end{align*}
We can rewrite the right-hand side and obtain
\begin{align*}
& \sum_{B_{1}\in \mathfrak{B}_1} \int_{  B_{s_i+1}\cap B_{1}} \prod_{n=1}^3 \Big( \sum_{T_{n}\in \mathcal{T}_n } \chi_{T_{n}}\Big)^{\frac{1}{2}}\lesssim N(B_{s_i+1}) \mathcal{K}(1,s_i)\fint _{P_{s_i}}\prod_{n=1}^3 \Big(\sum_{T_{n}\in \mathcal{T}_n(P_{s_i})} \chi_{\mathcal{T}_n(T_{n}, P_{s_i})}\Big)^{\frac{1}{2}}.
\end{align*}
Now, by construction, for $n=1,2$,
\begin{align*}
\fint _{P_{s_i}}&\Big(\sum_{T_{n}\in \mathcal{T}_n(P_{s_i})} \chi_{\mathcal{T}_n(T_{n}, P_{s_i})}\Big)^{\frac{1}{2}}\lesssim N(B_{s_i+1})^{-\frac{1}{2}}\beta_{n,2,s_i+1}^{\frac{1}{2}}\fint _{B_{s_i+1}}\Big(\sum_{T_{n}\in\mathcal{T}_n(B_{s_i+1})} \chi_{\mathcal{T}_n(T_{n},B_{s_i+1})}\Big)^{\frac{1}{2}}
\end{align*}
and for $n=3$,
\begin{align*}
\fint _{P_{s_i}}\Big(\sum_{T_{3}\in \mathcal{T}_3(P_{s_i})} \chi_{\mathcal{T}_3(T_{3}, P_{s_i})}\Big)^{\frac{1}{2}}&\sim |\mathcal{T}_3(P_{s_i})|^{\frac{1}{2}} |\mathcal{T}_3(B_{s_i+1})|^{-\frac{1}{2}} \fint _{B_{s_i+1}}\Big(\sum_{T_{3}\in\mathcal{T}_3(B_{s_i+1})} \chi_{\mathcal{T}_3(T_{3},B_{s_i+1})}\Big)^{\frac{1}{2}}\\
&\lesssim \mu_{4,s_{i}+1}^{-\frac{1}{2}}\mu_{5,s_{i}+1}^{-\frac{1}{2}}\beta_{3,2,s_i+1}^{\frac{1}{2}}\fint _{B_{s_i+1}}\Big(\sum_{T_{3}\in\mathcal{T}_3(B_{s_i+1})} \chi_{\mathcal{T}_3(T_{3},B_{s_i+1})}\Big)^{\frac{1}{2}},\end{align*}
where we used the fact that, by construction, $$|\mathcal{T}_3(B_{s_i+1})|= L(B_{s_i+1}\cap \mathfrak{A}_{1,s_i}) |\mathcal{T}_3(P_{s_i})|  \beta_{3,2,s_i+1}^{-1} \mu_{4,s_{i}+1}=\mu_{5,s_{i}+1}  |\mathcal{T}_3(P_{s_i})|  \beta_{3,2,s_i+1}^{-1} \mu_{4,s_{i}+1}.$$ Hence,
\begin{align*}
& \sum_{B_{1}\in \mathfrak{B}_1} \int_{  B_{s_i+1}\cap B_{1}} \prod_{n=1}^3 \Big( \sum_{T_{n}\in \mathcal{T}_n } \chi_{T_{n}}\Big)^{\frac{1}{2}}\\
&\lesssim \mu_{4,s_i+1}^{-\frac{1}{2}}\mu_{5,s_i+1}^{-\frac{1}{2}} \overline{\pmb{\beta}}_{2,s_i+1}^{\frac{1}{2}}\mathcal{K}(1,s_i) \prod_{n=1}^3\fint _{B_{s_i+1}}\Big(\sum_{T_{n}\in\mathcal{T}_n(B_{s_i+1})} \chi_{\mathcal{T}_n(T_{n},B_{s_i+1})}\Big)^{\frac{1}{2}}\\
&\sim  \mu_{4,s_{i}+1}^{-\frac{1}{2}}\mu_{5,s_{i}+1}^{-\frac{1}{2}} \overline{\pmb{\beta}}_{2,s_i+1}^{\frac{1}{2}}\mathcal{K}(1,s_i) \fint _{B_{s_i+1}}\prod_{n=1}^3\Big(\sum_{T_{n}\in\mathcal{T}_n(B_{s_i+1})} \chi_{\mathcal{T}_n(T_{n},B_{s_i+1})}\Big)^{\frac{1}{2}}.
\end{align*}
Now, if $s_{i}+1=S$, then, since $\mathcal{S}(\pmb{\beta},\pmb{\mu},S)=\mu_{4,S}^{-\frac{1}{2}}\mu_{5,S}^{-\frac{1}{2}} \overline{\pmb{\beta}}_{2,S}^{\frac{1}{2}}$ and $B_S$ is the ball of radius $R$, 
\begin{align*}
\fint _{B_{S}}\prod_{n=1}^3\Big(\sum_{T_{n}\in\mathcal{T}_n(B_{S})} \chi_{\mathcal{T}_n(T_{n},B_{S})}\Big)^{\frac{1}{2}}\sim  \prod_{n=1}^3|\mathcal{T}_n (Q_R)|^{\frac{1}{2}},
\end{align*}
and we are done. If $2\leq s_{i}+1\leq S-1$, we observe that
\begin{align*}
\fint _{B_{s_i+1}}&\prod_{n=1}^3\Big(\sum_{T_{n}\in\mathcal{T}_n(B_{s_i+1})} \chi_{\mathcal{T}_n(T_{n},B_{s_i+1})}\Big)^{\frac{1}{2}}\sim \fint _{B_{s_i+1}}\prod_{n=1}^3\Big(\sum_{T_{n}\in\mathcal{T}_n[M_{s_i+1} j+N_{s_i+1} t](B_{s_i+1},B_{s_i+1})} \chi_{T_{n}}\Big)^{\frac{1}{2}}.
\end{align*}
For any $P_{s_{i+1}}\in  \mathfrak{A}_{1,s_{i+1}}$, we have
\begin{align*}
& \sum_{B_{1}\in \mathfrak{B}_1} \int_{ P_{s_{i+1}}\cap B_{1}} \prod_{n=1}^3 \Big( \sum_{T_{n}\in \mathcal{T}_n } \chi_{T_{n}}\Big)^{\frac{1}{2}}\\
&\sim  \sum_{B_{s_i+1}\in  \mathfrak{A}_{2,s_i+1}}  \sum_{B_{1}\in \mathfrak{B}_1} \int_{P_{s_{i+1}}\cap B_{s_i+1}\cap B_{1}} \prod_{n=1}^3 \Big(\sum_{T_{n}\in \mathcal{T}_n } \chi_{T_{n}}\Big)^{\frac{1}{2}}\\
&\lesssim  \mu_{4,s_i+1}^{-\frac{1}{2}}\mu_{5,s_i+1}^{-\frac{1}{2}}\overline{ \pmb{\beta}}_{2,s_i+1}^{\frac{1}{2}}\mathcal{K}(1,s_i) |B_{s_i+1}|^{-1} \\
& \hspace{20mm}\sum_{B_{s_i+1}\in  \mathfrak{A}_{2,s_i+1}}\int _{P_{s_{i+1}}\cap B_{s_i+1}}\prod_{n=1}^3\Big(\sum_{T_{n}\in\mathcal{T}_n[M_{s_i+1} j+N_{s_i+1} t](B_{s_i+1},B_{s_i+1})} \chi_{T_{n}}\Big)^{\frac{1}{2}}\\
&\sim  \mu_{4,s_i+1}^{-\frac{1}{2}}\mu_{5,s_i+1}^{-\frac{1}{2}} \overline{\pmb{\beta}}_{2,s_i+1}^{\frac{1}{2}}\mathcal{K}(1,s_i)  |B_{s_i+1}|^{-1} \\
&\hspace{20mm} \sum_{B_{s_i+1}\in  \mathfrak{A}_{2,s_i+1}}\int _{P_{s_{i+1}}\cap B_{s_i+1}}\prod_{n=1}^3\Big(\sum_{T_{n}\in \mathcal{T}_n[M_{s_i+1} j+N_{s_i+1} t]( \mathfrak{A}_{2,s_i+1},P_{s_{i+1}})} \chi_{T_{n}}\Big)^{\frac{1}{2}}.
\end{align*}
By the definition of $\mathcal{K}(s_{i}+1,s_{i+1}) $,
\begin{align*}
 \sum_{B_{1}\in \mathfrak{B}_1} &\int_{ P_{s_{i+1}}\cap B_{1}} \prod_{n=1}^3 \Big( \sum_{T_{n}\in \mathcal{T}_n } \chi_{T_{n}}\Big)^{\frac{1}{2}}\\
\lesssim &\;\mu_{4,s_i+1}^{-\frac{1}{2}}\mu_{5,s_i+1}^{-\frac{1}{2}} \overline{\pmb{\beta}}_{2,s_i+1}^{\frac{1}{2}}\mathcal{K}(1,s_i)  |B_{s_i+1}|^{-1} \mathcal{K}(s_i+1,s_{i+1}) 2^{3(M_{s_i+1}j+N_{s_i+1} t)}   \\
&\prod_{n=1}^3| \mathcal{T}_n[M_{s_i+1} j+N_{s_i+1} t](\mathfrak{A}_{2,s_i+1},P_{s_{i+1}})|^{\frac{1}{2}}.\end{align*}
Now, applying Lemma \ref{L1}, we have
\begin{align*}
 \sum_{B_{1}\in \mathfrak{B}_1} &\int_{ P_{s_{i+1}}\cap B_{1}} \prod_{n=1}^3 \Big( \sum_{T_{n}\in \mathcal{T}_n } \chi_{T_{n}}\Big)^{\frac{1}{2}}\\
\lesssim &\;\mathcal{S}(\pmb{\beta},\pmb{\mu},s_{i}+1)\mathcal{K}(1,s_i)  |B_{s_i+1}|^{-1} 2^{3(M_{s_i+1}j+N_{s_i+1} t)}   \\
&\big(2^{2(1-F_{s_i+1})t}\big)^{\frac{1}{4}} \big(2^{(1-E_{s_i+1})j}\big)^{\frac{1}{2}}\prod_{n=1}^3|  \mathcal{T}_n[M_{s_i+1} j+N_{s_i+1} t](\mathfrak{A}_{2,s_i+1},P_{s_{i+1}})|^{\frac{1}{2}}.\end{align*}
A computation gives \begin{align*} |B_{s_i+1}|^{-1}  2^{3(M_{s_i+1}j+N_{s_i+1} t)}  &\sim   |B_{s_i+1}|^{-1} 2^{3((M_{s_i} +E_{s_i+1})j+(N_{s_i}+F_{s_i+1}) t)}  \\ &\sim 2^{-2(1-F_{s_i+1})t}2^{-(1-E_{s_i+1})j}\end{align*} and $$ \prod_{n=1}^3|\mathcal{T}_n[M_{s_i+1} j+N_{s_i+1} t](\mathfrak{A}_{2,s_i+1},P_{s_{i+1}})|^{\frac{1}{2}}\sim   2^{\frac{3}{2}(1-F_{s_i+1})t} 2^{\frac{1}{2}(1-E_{s_i+1})j}\prod_{n=1}^3|\mathcal{T}_n(P_{s_{i+1}})|^{\frac{1}{2}}.$$
Therefore,
\begin{align*}
& \sum_{B_{1}\in \mathfrak{B}_1} \int_{ P_{s_{i+1}}\cap B_{1}} \prod_{n=1}^3 \Big( \sum_{T_{n}\in \mathcal{T}_n } \chi_{T_{n}}\Big)^{\frac{1}{2}}\lesssim \mathcal{S}(\pmb{\beta},\pmb{\mu},s_{i}+1)\mathcal{K}(1,s_i)  \prod_{n=1}^3|\mathcal{T}_n(P_{s_{i+1}})|^{\frac{1}{2}},\end{align*}
which is what we wanted to prove.
\begin{flushright}
$\blacksquare$
\end{flushright}
\textbf{Proof of Lemma \ref{L1}.}

For the sake of notational compactness, we write $\mathcal{T}_n^*= \mathcal{T}_n[M_{s}j+N_{s}t]( \mathfrak{B}_{s}, P_{s'})$, and denote by $Q(s)$ the ball of radius $2^{M_{s}j+N_{s}t}$. We define $A_1 =\mu_{1,s}^{\frac{1}{4}}\big(2^{t(1-F_{s})}\big)^{\frac{1}{2}} $, $A_2=\mu_{2,s}^{\frac{1}{2}} \big(2^{j(1-E_{s})}\big)^{\frac{1}{2}}$ and $A_3=\mu_{3,s}$. 
We need to prove for $P_{s'}\in\mathfrak{A}_{1,s'}$,
\begin{align*}
 &\sum_{B_{s}\in \mathfrak{B}_s} \int_{ P_{s'}\cap B_{s}} \prod_{n=1}^3 \Big( \sum_{T_{n}\in \mathcal{T}_n^* } \chi_{T_{n}}\Big)^{\frac{1}{2}}\lesssim  |Q(s)| \overline{\pmb{\beta}}_{1,s}^{\frac{1}{2}}A_1 A_2A_3 \prod_{n=1}^3 |\mathcal{T}_n^*|^{\frac{1}{2}}.
\end{align*}
Now, as
\begin{align}\label{nok}
 \sum_{B_{s}\in \mathfrak{B}_s} \int_{ P_{s'}\cap B_{s}} \prod_{n=1}^3 \Big( \sum_{T_{n}\in \mathcal{T}_n^* } \chi_{T_{n}}\Big)^{\frac{1}{2}}\lesssim   |Q(s)|A_1^4 A_2^2 A_3 \Big(\prod_{n=1}^3 R_{n,1,s}^{\frac{1}{2}}\Big) |\mathfrak{A}_{1,s}\cap P_{s'}|,
\end{align}
it suffices to prove that \begin{align*}
 &  A_1^4 A_2^2  \Big(\prod_{n=1}^3 R_{n,1,s}^{\frac{1}{2}}\Big) 2^{-2t}2^{-j} \Big(\prod_{n=1}^3 R_{n,2,s}^{-\frac{1}{2}}\Big) \sum_{P\in \mathfrak{A}_{1,s}}\int_{ P_{s'}\cap P}\prod_{n=1}^3 \Big(\sum_{T_{n}\in \mathcal{T}_n[M_{s}j+N_{s}t](\mathfrak{A}_{1,s},P_{s'})}\chi_{T_{n}}\Big)^{\frac{1}{2}}\\
&\lesssim  |Q(s)| \overline{\pmb{\beta}}_{1,s}^{\frac{1}{2}}A_1 A_2 \prod_{n=1}^3 |\mathcal{T}_n^*|^{\frac{1}{2}}.
\end{align*}
It is easy to check that $$\prod_{n=1}^3|\mathcal{T}_n[M_{s}j+N_{s}t](\mathfrak{A}_{1,s}, P_{s'})|^{\frac{1}{2}}\lesssim 2^{\frac{3}{2}t-\frac{1}{2}r+\frac{1}{2}j}\prod_{n=1}^3|\mathcal{T}_n^*|^{\frac{1}{2}}.$$ This, together with the trilinear Kakeya estimate \eqref{GK1}, yields
\begin{align*}
 & A_1^4 A_2^2  \Big(\prod_{n=1}^3 R_{n,1,s}^{\frac{1}{2}}\Big) 2^{-2t}2^{-j} \Big(\prod_{n=1}^3 R_{n,2,s}^{-\frac{1}{2}}\Big) \sum_{P\in \mathfrak{A}_{1,s}}\int_{ P_{s'}\cap P}\prod_{n=1}^3 \Big(\sum_{T_{n}\in \mathcal{T}_n[M_{s}j+N_{s}t](\mathfrak{A}_{1,s},P_{s'})}\chi_{T_{n}}\Big)^{\frac{1}{2}}\\
& \lesssim  A_1^4 A_2^2 \Big(\prod_{n=1}^3 R_{n,1,s}^{\frac{1}{2}}\Big) 2^{-2t}2^{-j} \Big(\prod_{n=1}^3 R_{n,2,s}^{-\frac{1}{2}}\Big)  2^{\frac{1}{2}(j+t+r)} |Q(s)| \prod_{n=1}^3|\mathcal{T}_n[M_{s}j+N_{s}t](\mathfrak{A}_{1,s}, P_{s'})|^{\frac{1}{2}}\\
& = A_1^3 A_2 \Big(\prod_{n=1}^3 R_{n,1,s}^{\frac{1}{2}}\Big)A_1 A_2\Big(\prod_{n=1}^3 R_{n,2,s}^{-\frac{1}{2}}\Big)  |Q(s)| \prod_{n=1}^3|\mathcal{T}_n^*|^{\frac{1}{2}}.
\end{align*}

We now present the following key lemma:
\begin{Lemma}\label{cord}
\begin{align*}
A_1^3 A_2 \Big(\prod_{n=1}^3 R_{n,1,s}^{\frac{1}{2}}\Big)\lesssim  \overline{\pmb{\beta}}_{1,s}^{\frac{1}{2}}\Big(\prod_{n=1}^3 R_{n,2,s}^{\frac{1}{2}}\Big) .
\end{align*}
\end{Lemma}

\textbf{Proof of Lemma \ref{cord}.}

We have for any $P\in \mathfrak{A}_{1,s}$,
\begin{align}\label{eq1}\nonumber
 |Q(s)|A_1^4 A_2^2 A_3  R_{1,1,s} R_{3,1,s}&\sim \sum_{B_{s}\in \mathfrak{B}_s} \int_{P\cap B_s} \Big(\sum_{T_{1}\in \mathcal{T}_1^*} \chi_{T_{1}}\Big)\Big(\sum_{T_{3}\in \mathcal{T}_3^*} \chi_{T_{3}}\Big)\\\nonumber
&\sim 2^{-r}A_3\beta_{1,1,s} \beta_{3,1,s} \int_{P}  \Big(\sum_{T_{1}\in \mathcal{T}_1^*} \chi_{T_{1}}\Big)\Big(\sum_{T_{3}\in \mathcal{T}_3^*} \chi_{T_{3}}\Big)\\
&\lesssim  |Q(s)| A_3\beta_{1,1,s} \beta_{3,1,s}R_{1,2,s} R_{3,2,s}.
\end{align}
In the same way, we obtain
\begin{align}\label{eq2}
A_1^4 A_2^2 R_{2,1,s} R_{3,1,s}\lesssim  2^r\beta_{1,2,s}\beta_{3,2,s}R_{2,2,s} R_{3,2,s}.
\end{align}
In the other case, we have
\begin{align}\label{eq3}\nonumber
 |Q(s)|A_1^4 A_2^2  A_3 R_{1,1,s} R_{2,1,s}&\sim\sum_{B_{s}\in \mathfrak{B}_s} \int_{P\cap B_s} \Big(\sum_{T_{1}\in \mathcal{T}_1^*} \chi_{T_{1}}\Big)\Big(\sum_{T_{2}\in \mathcal{T}_2^*} \chi_{T_{2}}\Big)\\\nonumber
&\sim 2^{-j}A_2^2 A_3\beta_{1,1,s} \beta_{2,1,s}  \int_{P}  \Big(\sum_{T_{1}\in \mathcal{T}_1^*} \chi_{T_{1}}\Big)\Big(\sum_{T_{2}\in \mathcal{T}_2^*} \chi_{T_{2}}\Big)\\
&\lesssim  |Q(s)|A_2^2 A_3 \beta_{1,1,s} \beta_{2,1,s} R_{1,2,s} R_{2,2,s}.
\end{align}
Combining \eqref{eq1}, \eqref{eq2} and \eqref{eq3},
\begin{align*}
A_1^4 A_2^{\frac{4}{3}} \Big(\prod_{n=1}^3 R_{n,1,s}^{\frac{2}{3}} \Big)\lesssim  \overline{\pmb{\beta}}_{1,s}^{\frac{2}{3}}\Big(\prod_{n=1}^3R_{n,2,s}^{\frac{2}{3}}\Big) .
\end{align*}
The result follows by raising both sides to the power $\frac{3}{4}$.
\begin{flushright}
$\blacksquare$
\end{flushright}

Using the previous lemma, we conclude
\begin{align*}
&A_1^3 A_2  \Big(\prod_{n=1}^3 R_{n,1,s}^{\frac{1}{2}}\Big)A_1 A_2 \Big(\prod_{n=1}^3 R_{n,2,s}^{-\frac{1}{2}}\Big)  |Q(s)| \prod_{n=1}^3|\mathcal{T}_n^*|^{\frac{1}{2}}\\& \lesssim  2^{\frac{1}{2}r}  \overline{\pmb{\beta}}_{1,s}^{\frac{1}{2}}A_1 A_2  |Q(s)| \prod_{n=1}^3|\mathcal{T}_n^*|^{\frac{1}{2}},
\end{align*}
as desired.

\begin{flushright}
$\blacksquare$
\end{flushright}

\begin{Remark}\label{noKak}
We observe that if $s'=s$, then \eqref{nok}, together with Lemma \ref{cord}, implies Lemma \ref{L2}. In particular, the proof of Lemma \ref{L2} for the case $s'=s$ does not rely on the trilinear Kakeya estimate \eqref{GK1}.
\end{Remark}

\subsection{An alternative proof of the trilinear Kakeya estimate with sharp dependence on transversality}
In this section, we reprove Guth's trilinear Kakeya estimate with sharp dependence on transversality. 
\begin{Theorem}\label{MKe}
Fix $r,j,t,w,m$ and let $\mathcal{T}=\{\mathcal{T}_n\}_n$ be a collection of tubes of $(r,j,t,w,m) $-type. Then given any $\epsilon>0$ there exist a constants $C_\epsilon<\infty$ such that
\begin{align*}
\int_{Q_R} \prod_{n=1}^3 \Big( \sum_{T_{n}\in \mathcal{T}_n } \chi_{T_{n}}\Big)^{\frac{1}{2}}\leq C_\epsilon R^{\epsilon} 2^{\frac{1}{2}(j+t+r)}  \prod_{n=1}^3 |\mathcal{T}_n|^{\frac{1}{2}},
\end{align*}
for all $R>0$.
\end{Theorem}

Our proof is slightly weaker than Guth's \cite{G}, as it requires the factor $R^{\epsilon}$, but it is simpler and follows a philosophy similar to that of \cite{gu2} and the arguments from the previous section. We define:

\begin{Definition}
Fix $s,r, j,t,w,m$ and let $\mathcal{T}=\{\mathcal{T}_n\}_n$ be a collection of tubes of $(r,j,t,w,m) $-type. We denote by $\mathcal{M}(s)$ the smallest constant $C$ such that 
\begin{align*}
\int_{P_s} \prod_{n=1}^3 \Big( \sum_{T_{n}\in \mathcal{T}_n } \chi_{T_{n}}\Big)^{\frac{1}{2}}\leq C 2^{\frac{1}{2}(j+t+r)} \prod_{n=1}^3 |\mathcal{T}_n(P_{s})|^{\frac{1}{2}},
\end{align*}
for every $P_{s}\in  \mathcal{A}(Q_R, \mathcal{P}(j,t)[s])$. 
\end{Definition}

As in \cite{gu2}, we can use the Loomis-Whitney inequality \cite{lw} to obtain, after an affine transformation, the following:
\begin{Lemma}\label{WLK} (Loomis-Whitney)
For any three directions $v_1, v_2, v_3$ such that $|v_1\wedge v_2\wedge v_3|=\theta>0$, and collections of tubes $\{\mathcal{T}_n\}_{n=1}^3$ such that each tube $T_n\in\mathcal{T}_n$ satisfies $d(T_n)=v_n$, we have
\begin{align}\label{LW}
\int_{\R^d}\prod_{n=1}^3 \Big( \sum_{T_{n}\in \mathcal{T}_n } \chi_{T_{n}}\Big)^{\frac{1}{2}}\lesssim \theta^{-\frac{1}{2}}  \prod_{n=1}^3 |\mathcal{T}_n|^{\frac{1}{2}}.
\end{align}
 \end{Lemma}

Let $\{\mathcal{T}_n\}_{n=1}^3$ be $3$ families of tubes of $(r,j,t,w,m) $-type. Let  $S_1, S_2, S_3$ be the subsets of the paraboloid $\{(\xi,\frac{1}{2}|\xi|^2): \; |\xi|\leq 2\}$ such that for every $T_n\in \mathcal{T}_n$, we have $d(T_n)=n(\xi_n)$ for some $\xi_n\in S_n$. We can decompose 
$$\mathcal{T}_n=\bigcup_{i=1}^N \mathcal{T}_{n,i} $$
where $N\sim 2^{2\lambda}$, and if $T_n\in \mathcal{T}_{n,i} $, then $d(T_n)=n(\xi_n)$, for some $\xi_n\in S_{n,i}$, with
\begin{align*}
&S_{1,i}\subset \{(\xi,\tfrac{1}{2}|\xi|^2):\;\xi\in\tau^{j+\lambda}_{k_i} \cap \mathfrak{t}^{j+t+\lambda}_{w,m_i}\},\\
&S_{2,i}\subset \{(\xi,\tfrac{1}{2}|\xi|^2):\;\xi\in\tau^{j+\lambda}_{k_i'} \cap \mathfrak{t}^{j+t+\lambda}_{w,m_i}\},\\
&S_{3,i}\subset\{(\xi,\tfrac{1}{2}|\xi|^2):\;\xi\in\tau^{r+\lambda}_{k_i''} \cap \mathfrak{t}^{r+t+\lambda}_{w',m_i},\},
\end{align*}
for some indices $k_i,k_i',k_i'',m_i$. 

Taking a parameter $\lambda$ such that $\lambda\sim \frac{1}{\epsilon} \log_2 C$, for some constant $C$ independent on $R$ and $\theta$, it suffices to prove Theorem \ref{MKe} for each collection of tubes $\{\mathcal{T}_{n,i}\}_n$.

\begin{Remark}\label{re} Let $P_{s+\lambda}\in  \mathcal{A}(Q_R, \mathcal{P}(j,t)[s+\lambda])$. We can see that for every $P_s\in \mathcal{A}(Q_R, \mathcal{P}(j,t)[s])\cap P_{s+\lambda}$, $T_n[s], T_{n}[s]'\in \mathcal{T}_{n,i}[s](P_s)$, we have \begin{align*} \mathcal{T}_{n,i}[s](T_{n}[s],P_s)\cap P_{s+\lambda}\subset \bigcup_{a\in A}\big(\mathcal{T}_{n,i}[s](T_{n}[s]',P_s)+a\big)\cap P_{s+\lambda},\end{align*}
for some set of points $A$ such that $|A|\lesssim 1$.

\end{Remark}

Theorem \ref{MKe} follows from the following lemmas:

\begin{Lemma}\label{LK1}For every $s\geq 0$, $$M(s+\lambda)\lesssim M(s).$$\end{Lemma}\label{L2C}\begin{Lemma}\label{LK2}$$M(0)\lesssim 1.$$\end{Lemma}

\textbf{Proof of Theorem \ref{MKe}.}

It is enough to prove that $\mathcal{M}(\log_2 R)\leq C_\epsilon R^\epsilon$. By Lemma \ref{LK1} applied $\frac{\log_2 R}{\lambda}$ times:
\begin{align*}
\mathcal{M}(\log_2 R)&\leq C \mathcal{M}(\log_2 R-\lambda)\\
&\leq C^2 \mathcal{M}(\log_2 R-2\lambda)\\
&\leq \cdots \leq C^{\frac{\log_2 R}{\lambda}} \mathcal{M}(0)\\
&= R^{\frac{\log_2 C}{\lambda}} \mathcal{M}(0)\\
\end{align*}
Using Lemma \ref{LK2} and that $\lambda\sim \frac{1}{\epsilon} \log_2 C$, we conclude the result.
\begin{flushright}
$\blacksquare$
\end{flushright}

\textbf{Proof of Lemma \ref{LK1}.}
Let $P_{s+\lambda}\in  \mathcal{A}(Q_R, \mathcal{P}(j,t)[s+\lambda])$, we have by definition,
\begin{align*}
\int_{P_{s+\lambda}} \prod_{n=1}^3 \Big( \sum_{T_{n}\in \mathcal{T}_n } \chi_{T_{n}}\Big)^{\frac{1}{2}}&\leq \sum_{P_s\in  \mathcal{A}(Q_R, \mathcal{P}(j,t)[s])}\int_{P_{s+\lambda}\cap P_{s}} \prod_{n=1}^3 \Big( \sum_{T_{n}\in \mathcal{T}_n } \chi_{T_{n}}\Big)^{\frac{1}{2}}\\
&\lesssim \sum_{P_s\in  \mathcal{A}(Q_R, \mathcal{P}(j,t)[s])} M(s) 2^{\frac{1}{2}(t+j+r)}  \prod_{n=1}^3 |\mathcal{T}_n(P_{s}\cap P_{s+\lambda})|^{\frac{1}{2}}.
\end{align*}
As in the previous section, we observe that  \begin{align}\label{Tubr}\prod_{n=1}^3 |\mathcal{T}_n[s](\mathcal{A}(Q_R, \mathcal{P}(j,t)[s]),P_{s+\lambda})|^{\frac{1}{2}}\sim 2^{\frac{1}{2}j-\frac{1}{2}r+\frac{3}{2}t}\prod_{n=1}^3 |\mathcal{T}_n(P_{s+\lambda})|^{\frac{1}{2}}.\end{align}
Thanks to the equation
\begin{align*}
|\mathcal{T}_n(P_{s}\cap P_{s+\lambda})|^{\frac{1}{2}}\sim \fint_{P_s}\prod_{n=1}^3 \Big(\sum_{T_n\in \mathcal{T}_n[s](\mathcal{A}(Q_R, \mathcal{P}(j,t)[s]),P_{s+\lambda})} \chi_{\mathcal{T}_n}\Big)^{\frac{1}{2}},
\end{align*}
we can rewrite the previous inequality as
\begin{align*}
\int_{P_{s+\lambda}} &\prod_{n=1}^3 \Big( \sum_{T_{n}\in \mathcal{T}_n } \chi_{T_{n}}\Big)^{\frac{1}{2}}\\
&\lesssim \sum_{P_s\in\mathcal{A}(Q_R, \mathcal{P}(j,t)[s])}M(s) 2^{\frac{1}{2}(t+j+r)}  \fint_{P_s}\prod_{n=1}^3 \Big(\sum_{T_n\in\mathcal{T}_n[s](\mathcal{A}(Q_R, \mathcal{P}(j,t)[s]),P_{s+\lambda})} \chi_{\mathcal{T}_n}\Big)^{\frac{1}{2}}\\
&\sim M(s)  2^{\frac{1}{2}(t+j+r)}|P_s|^{-1}  \int_{P_{s+\lambda}}\prod_{n=1}^3 \Big(\sum_{T_n\in\mathcal{T}_n[s](\mathcal{A}(Q_R, \mathcal{P}(j,t)[s]),P_{s+\lambda})} \chi_{\mathcal{T}_n}\Big)^{\frac{1}{2}}.
\end{align*}
By Remark \ref{re}, we can apply Lemma \ref{WLK}, which together with \eqref{Tubr} yields
\begin{align*}
\int_{P_{s+\lambda}} &\prod_{n=1}^3 \Big( \sum_{T_{n}\in \mathcal{T}_n } \chi_{T_{n}}\Big)^{\frac{1}{2}}\\
&\lesssim  M(s) 2^{\frac{1}{2}(t+j+r)}|P_s|^{-1}  2^{3s}2^{\frac{1}{2}(t+j+r)} \prod_{n=1}^3 |\mathcal{T}_n[s](\mathcal{A}(Q_R, \mathcal{P}(j,t)[s]),P_{s+\lambda})|^{\frac{1}{2}}\\
&\lesssim M(s)2^{\frac{1}{2}(t+j+r)}|P_s|^{-1} 2^{3s} 2^{\frac{1}{2}(t+j+r)}2^{\frac{1}{2}j-\frac{1}{2}r+\frac{3}{2}t}\prod_{n=1}^3 |\mathcal{T}_n(P_{s+\lambda})|^{\frac{1}{2}}.
\end{align*}
As we have $|P_s|\sim 2^{3s} 2^{2t+j}$, we conclude that
\begin{align*}
\int_{P_{s+\lambda}} \prod_{n=1}^3 \Big( \sum_{T_{n}\in \mathcal{T}_n } \chi_{T_{n}}\Big)^{\frac{1}{2}}\lesssim M(s) 2^{\frac{1}{2}(t+j+r)} \prod_{n=1}^3 |\mathcal{T}_n(P_{s+\lambda})|^{\frac{1}{2}},
\end{align*}
as desired.

\begin{flushright}
$\blacksquare$
\end{flushright}

\textbf{Proof of Lemma \ref{LK2}.}
It follows directly from Lemma \ref{L1} and Remark \ref{noKak}.
\begin{flushright}
$\blacksquare$
\end{flushright}

\end{document}